\documentclass[a4paper,12pt,reqno]{amsart}
\usepackage{amssymb}
\usepackage{amsmath}

\def\bu{{\hskip-1pt}_\star{\hskip-1pt}}
\def\i{\,\lrcorner\,}
\def\a{\alpha}

\def\ni{\noindent}
\def\vs{\vskip .6cm}

\def\la{\langle}
\def\ra{\rangle}
\def\.{\cdot}
\def\O{\Omega}
\def\n{\nabla}
\def\l{\lambda}
\def\s{\sigma}

\def\beq{\begin{equation}}
\def\eeq{\end{equation}}
\def\bea{\begin{eqnarray*}}
\def\eea{\end{eqnarray*}}
\def\beaa{\begin{eqnarray}}
\def\eeaa{\end{eqnarray}}
\def\ba{\begin{array}}
\def\ea{\end{array}}
\def\x{\times}
\def\f{\varphi}

\def\o{\omega}

\def\L{\Lambda}

\def\bp{\begin{proof}}
\def\r{\end{proof}}

\def\D{\mathfrak D}

\def \RM{\mathbb{R}}

\def\End{{\rm End}}

\def\d{{\delta}}

\def\Ric{\mathrm{Ric}}

\def\id{\mathrm{id}}
\def\be{\begin{equation}}

\def\ee{\end{equation}}
\def\Lie{{\mathcal L\,}}

\def\tr{\mathrm{tr}}

\def\RP{\mathbb{R}\mathrm{P}}

\def\Sym{\mathrm{Sym}}

\def\R{\mathbb{R}}

\def\SU{\mathrm{SU}}
\def\SO{\mathrm{SO}}
\def\Sp{\mathrm{Sp}}
\def\U{\mathrm{U}}

\def\M{\mathfrak{M}}

\def\scal{\mathrm{scal}}

\def\psp{\psi ^+}
\def\psm{\psi ^-}
\def\dpsp{\dot\psi ^+}
\def\dpsm{\dot\psi ^-}

\newtheorem{ede}{Definition}[section]

\newtheorem{ath}[ede]{Theorem}

\newtheorem{elem}[ede]{Lemma}

\title{Deformations of Nearly K\"ahler Structures}

\author{Andrei Moroianu, Paul-Andi Nagy and Uwe Semmelmann}

\address{Andrei Moroianu \\ CMAT\\ {\'E}cole Polytechnique \\ UMR 7640 du CNRS
\\ 91128 Palaiseau \\ France}
\email{am@math.polytechnique.fr}

\address{Paul-Andi Nagy \\
Department of Mathematics \\
University of Auckland \\
Private Bag 92019\\ Auckland\\ New Zealand}
\email{nagy@maths.auckland.ac.nz}

\address{Uwe Semmelmann\\ Mathematisches Institut, Universit{\"a}t zu
K{\"o}ln\\
Weyertal 86-90 D-50931 K{\"o}ln, Germany}
\email{uwe.semmelmann@math.uni-koeln.de}

\begin{document}

\begin{abstract}
We study the space of nearly K\"ahler structures on compact
$6$-dimensional manifolds. In particular, we prove that the 
space of infinitesimal deformations of a strictly nearly K\"ahler
structure (with scalar curvature $\scal$) modulo the
group of diffeomorphisms, is isomorphic to the space of primitive co-closed
$(1,1)$-eigenforms of the Laplace operator for the eigenvalue $2\scal/5$.

\date\today
\vs

\noindent
2000 {\it Mathematics Subject Classification}: Primary 58E30, 53C10, 53C15.

\medskip
\noindent{\it Keywords:} deformations, $\SU_3$ structures, nearly
K\"ahler manifolds, Gray manifolds
\end{abstract}

\maketitle

\section{Introduction}

A nearly K\"ahler manifold is an almost Hermitian
manifold $(M,g,J)$ with the property that $(\n_X J)X=0$ for all tangent vectors
$X$, where $\n$ denotes the Levi-Civita connection of $g$. A nearly
K\"ahler manifold is called {\em strictly} nearly K\"ahler if $(\n_X
J)$ is non-zero for every tangent vector $X$. Besides
K\"ahler manifolds, there are  
two main families of examples of compact nearly K\"ahler  manifolds: Naturally
reductive $3$-symmetric spaces, which are classified by A. Gray and
J. Wolf \cite{gw} and twistor spaces over compact
quaternion-K\"ahler manifolds with positive scalar curvature, endowed
with the non-integrable canonical almost complex structure (cf. \cite{na}).

A nearly K\"ahler manifold of dimension $4$ is
automatically a K\"ahler surface, and the only known examples of
non-K\"ahler compact nearly K\"ahler manifolds in dimension $6$ are
the $3$-symmetric spaces $G_2/\SU_3$, $\SU_3/S^1\times S^1$,
$\Sp_2/S^1\x \Sp_1$ and $\Sp_1\x\Sp_1\x\Sp_1/\Sp_1$. Moreover,
J.-B. Butruille has recently shown in \cite{jbb}
that there are no other homogeneous examples
in dimension $6$.

On the other hand, using previous results of R. Cleyton and
A. Swann on $G$-structures with skew-symmetric intrinsic torsion, the
second-named author has proved in \cite{na} that every compact simply
connected nearly K\"ahler  manifold $M$ is isometric to a Riemannian product
$M_1\times\ldots\times M_k$, such that for each $i$, $M_i$ is a nearly
K\"ahler manifold belonging to the following list:
 K\"ahler manifolds,
naturally reductive $3$-symmetric spaces,
twistor spaces over compact quaternion-K\"ahler manifolds
  with positive scalar curvature, and $6$-dimensional nearly K\"ahler
  manifolds. 

It is thus natural to concentrate on the $6$-dimensional case, all
the more that in this dimension, non-K\"ahler nearly K\"ahler  manifolds have
several interesting
features: they carry real Killing spinors (and thus are automatically
Einstein with positive scalar curvature) and they are defined by a
$\SU_3$ structure whose intrinsic
torsion is skew-symmetric. These manifolds were intensively studied
by A. Gray in the 70's, thus motivating the following

\begin{ede} A compact strictly nearly K\"ahler  manifold of
  dimension $6$ is called a {\em Gray manifold}.
\end{ede}

The main goal of this paper is to study the deformation problem for Gray
manifolds. Notice that we consider simultaneous deformations of the
metric and of the 
almost complex structure. Indeed M.~Verbitsky
proved in \cite{mv} that on a 6-dimensional almost
complex manifold there is up to constant rescaling at most one
strictly nearly K\"ahler metric. Conversely it is well known
(c.f.~\cite{bfgk} or Section 4 below) that on a manifold $(M^6,g)$
which is not locally isometric to the standard sphere, there is at
most one compatible almost complex
structure $J$ such that $(M,g,J)$ is nearly K\"ahler.

We start by studying deformations of $\SU_3$ structures,
then use the characterization of Gray manifolds as $\SU_3$ structures
satisfying a certain exterior differential system in order to compute the space of
infinitesimal deformations of a given Gray structure modulo
diffeomorphisms. In particular, we prove that this space is isomorphic
to some eigenspace of the Laplace operator acting on $2$-forms (see
Theorem \ref{infd} for a precise statement).

\section{Algebraic preliminaries}

Let $V$ denote the standard $6$-dimensional $\SU_3$ representation
space, which comes equipped with the Euclidean product $g\in \Sym V^*$, the
complex structure $J\in \End(V)$, the fundamental $2$-form
$\o(\.,\.)=g(J\.,\.)\in\L^2 V^*$, and the complex volume element
$\psp+i\psm\in \L^{(3,0)}V^*$.

These objects satisfy the compatibility relations
\beq
\omega \wedge \psi^\pm = 0,\qquad
\psi^+ \wedge \psi^- = \frac23 \omega^3=4dv,
\eeq
where $dv$ denotes the volume form of the metric $g$.
It is easy to check that $\psp$ and $\psm$ are related by
\beq\label{psm}
\psm(X,Y,Z):=-\psp(JX,Y,Z).\eeq

We identify elements of $V$ and $V^*$ using the isomorphism induced by
$g$. For any orthonormal basis $\{e_i\}$ of $V$ adapted to $J$
({\em i.e.} $J(e_{2i-1})=e_{2i}$) we have:
$$\o=e ^{12}+e ^{34}+e ^{56},$$
$$\psp=e ^{135} - e ^{146} - e ^{236} - e ^{245}\ ,\qquad
\psm=e ^{136} +e ^{145} + e ^{235} - e ^{246}.$$
The following formulas are straightforward (it is enough to check them
for $X=e_1$ and use the transitivity of the $\SU_3$ action on spheres):
\beq\label{f1}\psp\wedge(X\i\psp)=X\wedge\o ^2,\qquad
\psp\wedge(X\i\psm)=-JX\wedge\o ^2,\eeq 
\beq\label{f3}\psm\wedge(X\i\psp)=JX\wedge\o ^2,\qquad
\psm\wedge(X\i\psm)=X\wedge\o ^2.\eeq 

Let $\L:\L^pV\to \L^{p-2}V$ denote the metric adjoint of the wedge
product with $\o$, $\Lambda =
\frac12\sum_i Je_i \lrcorner e_i \lrcorner$. It is easy to check
that 
\beq\label{lam}\L(X\i\psi ^\pm)=0,
\quad\L(X\wedge\psi ^\pm)=JX\i\psi ^\pm,\quad\forall X\in V \eeq 
and 
\beq\label{lam1}\L(\tau\wedge\o)=\o\wedge\L\tau+(3-p)\tau,\quad
\forall\tau\in\L^pV. \eeq

We next describe the decomposition into irreducible summands of
$\L^2V$ and $\L^3V$. We use the notation
$\Lambda^{(p,q)+(q,p)}V$ for the projection of $\Lambda^{(p,q)}V$ onto
the real space $\Lambda^{p+q}V$.
\beq
\Lambda^2V=(\Lambda^{(1,1)}_0V\oplus \R\omega) \oplus 
\Lambda^{(2,0)+(0,2)}V \ ,
\eeq
where the first two summands consist of $J$-invariant and the
last of $J$-anti-invariant forms. Here $\Lambda^{(1,1)}_0V$ is the
space of primitive $(1,1)$-forms, {\em i.e.} the kernel of the contraction
map $\Lambda$. The map 
\beq\label{ii}\xi \mapsto
\xi \lrcorner \psi^+\eeq 
defines an isomorphism of the second summand
$\Lambda^{(2,0)+(0,2)}V$ with $V$.  For
3-forms we 
have the irreducible decomposition
\beq\label{lambda3}
\Lambda^3V = (\Lambda^1V \wedge \omega) \oplus 
\Lambda^{(3,0)+(0,3)}V \oplus  \Lambda^{(2,1)+(1,2)}_0V.
\eeq
The second summand $\Lambda^{(3,0)+(0,3)}V$ is 2-dimensional and spanned
by the forms $\psi^\pm$. The third summand  $\Lambda^{(2,1)+(1,2)}_0V$
is 12-dimensional and can be identified with the space of symmetric
endomorphisms of $V$ anti-commuting with $J$. We note that because
of the Schur lemma, the map given by taking the wedge product with
$\omega$ vanishes on the last two summands.

An endomorphism $A$ of $V$ (not necessarily skew-symmetric) acts on
$p$-forms by the formula
\beq\label{np}(A\bu u)(X_1,\ldots,X_p):=-\sum_{i=1}^p
u(X_1,\ldots,A(X_i),\ldots,X_p).\eeq
A more convenient way of writing this action is
\beq\label{a}A\bu u=-\sum_{i=1}^p A ^*(e_i)\wedge e_i\i u,\eeq
where $A ^*$ denotes the metric adjoint of $A$. Taking $A=J$ we
obtain the form spaces $\Lambda^{(p,q)+(q,p)}V$ as eigenspaces
of the $J$-action for the eigenvalues $-(p-q)^2$. In particular,
we see that $J\bu \varphi = 0$ for any $J$-invariant 2-form  $\varphi \in
\Lambda^{(1,1)}V$.

Let $\Sym^-V$ denote the space of symmetric
endomorphisms anti-commuting with $J$. This space is clearly
invariant by composition with $J$. The map
$S \mapsto S\bu \psi^+$, with $S \in \Sym^-V$, defines an isomorphism
of $\SU(3)$-representations
$$
\Sym^-V\cong
\Lambda^{(2,1)+(1,2)}_0V,
$$
showing in particular that $\Sym^-V$ is irreducible. Taking (\ref{psm})
into account, we remark that for $S\in \Sym^-V$ we have
\beq\label{psps}S\bu \psp=(JS)\bu \psm.\eeq
Notice that $\tr(S)=0$ for all $S\in \Sym^-V$. 

Let $h$ be any skew-symmetric endomorphism anti-commuting with $J$. Then
the map $h \mapsto g(h\cdot,\cdot)$ identifies the space of
skew-symmetric endomorphism anti-commuting with $J$ with
$\Lambda^{(2,0)+(0,2)}V$. Using the isomorphism (\ref{ii}) we can state this as
\begin{elem}\label{ss} An endomorphism $F$ of $V$ anti-commuting with
  $J$ can be written in a unique way
$$F=S+\psp_\xi$$
for some $S\in \Sym^-V$ and $\xi\in V$, where $\psp_\xi$ denotes the
skew-symmetric 
endomorphism of $V$ defined by $g(\psp_\xi \.,\.)=\psp(\xi,\.,\.)$.
\end{elem}

Corresponding to the decomposition of $\Lambda^3V$ given in (\ref{lambda3}),
we have

\begin{elem}\label{ps} An exterior $3$-form $u\in\L^3V$ can be written in a
  unique way 
$$u=\a\wedge\o+\l \psp+\mu\psm+S\bu \psp,$$
for some $\a\in V$, $\l,\mu\in \RM$ and $S\in \Sym^-V$. Its contraction
with $\o$ satisfies 
\beq\label{cont}\L u=2\a.\eeq
\end{elem}
\bp The contraction map $\L$ obviously vanishes on
$\Lambda^{(3,0)+(0,3)}V\oplus \Lambda^{(2,1)+(1,2)}_0 V$, so by
(\ref{lam1}) we have $\L u=\L(\a\wedge\o)=2\a$.
\r

The space of symmetric endomorphisms commuting with $J$ is identified
with $\Lambda^{(1,1)}V$ via the map 
$h \mapsto \varphi(\cdot,\cdot) := g(hJ\cdot, \cdot)$, which in
particular maps the identity of $V$ to the fundamental form $\omega$.
If $\varphi$ is a $(1,1)$-form with 
corresponding symmetric endomorphism $h$, then
\beq
h = h_0 + \frac16\tr(h)\id,
\eeq
where $h_0$ denotes the trace-free part of $h$.
As a consequence of this formula and of the Schur Lemma we find that
\beq\label{hps}h\bu \psp=\frac16\tr(h)\id\bu\psp=-\frac12 \tr(h)\psp
\eeq
for all symmetric endomorphisms $h$ commuting with $J$.

In the remaining part of this section we want to recall several
properties and formulas related to the Hodge $*$-operator, which we
will use in later computations.
We consider the scalar product $\la\.,\.\ra$ on $\L^kV$ characterized
by the fact that the basis 
$$\{e_{i_1}\wedge\ldots\wedge e_{i_k}\ |\ 1\le i_1<\ldots<i_k\le 6\}$$
is orthonormal.
With respect to this scalar product, the interior and exterior products
are adjoint operators:
\begin{equation}\label{81}\la X\i\o,\tau\ra=\la\o,
  X\wedge\tau\ra,\qquad\forall\ X\in V,
\ \o\in\L^kV,\
 \tau\in\L^{k-1}V.
\end{equation}
We define the Hodge *-operator $*:\L^kV\to\L^{6-k}V$ by
$$\o\wedge *\tau:=\la\o,\tau\ra dv,\qquad\forall\ \o,\tau\in\L^kV,$$
where $dv=1/6 \omega^3$ denotes the volume form ($dv=e ^{123456}$ in
our notations). 
It is well-known  that the following relations are satisfied:
\beq\label{star}
*\o=\frac12\o ^2,\quad\la
*\o,*\tau\ra=\la\o,\tau\ra,\quad 
*^2=(-1)^{k} \quad \hbox{on}\  \L^kV
\eeq

From the expression of $\psp$ and $\psm$ in any orthonormal basis
$\{e_i\}$ we see that
$*\psp=\psm$ and $*\psm=-\psp$. For later use we compute the Hodge
operator on $\L^{(2,1)+(1,2)}_0V$, too. Let $S\in \Sym^-V$ and $\a\in
\L^{(2,1)+(1,2)}_0V$. We have $\a\wedge\psm=0$, whence
\bea \la\a,S\bu \psp\ra dv&\stackrel{(\ref{a})}{=}&-\sum_i \la
\a,S(e_i)\wedge e_i\i\psp\ra dv= -\sum_i \la e_i\wedge S(e_i)\i
\a,\psp\ra dv\\
&=& (S\bu \a)\wedge *\psp=  (S\bu \a)\wedge
\psm=S\bu (\a\wedge\psm)-\a\wedge(S\bu \psm)\\
&=&-\a\wedge(S\bu \psm)=\la\a,*(S\bu \psm)\ra dv \eea (we used the fact
that $*^2=-1$ on $3$-forms in order to get the last equality). This
shows that
\beq\label{sp}*(S\bu \psm)=S\bu \psp\qquad\hbox{and}\qquad*(S\bu \psp)=-S\bu \psm.\eeq
There are two other formulas which we will  use later. Let $\varphi_0$
be a primitive $(1,1)$-form, $\xi \in V$ and $\alpha \in
\Lambda^pV$ then 
\beq\label{sp1} *(\varphi_0 \wedge \omega) =
-\varphi_0,\qquad *(\xi \wedge\alpha) = (-1)^p\xi
\lrcorner *\alpha \ . \eeq

\vs

\section{Deformations of $\SU_3$ structures}

Let $M$ be  a smooth $6$-dimensional manifold. 

\begin{ede}
A $\SU_3$ structure on $M$ is a
reduction of the frame bundle of $M$ to 
$\SU_3$. It consists of a $5$-tuple
$(g,J,\o,\psp,\psm)$ where $g$ is a Riemannian metric, $J$ is a
compatible almost complex structure, $\o$ is the corresponding
fundamental $2$-form $\o(\.,\.)=g(J\.,\.)$ and $\psp+i\psm$ is a
complex volume form of type $(3,0)$. 
\end{ede}

If $M$ carries a $\SU_3$ structure, each tangent space $T_xM$ has a
$\SU_3$ representation isomorphic to the standard one, so all
algebraic results of the previous section transpose {\em verbatim} 
to global results on $M$. In the remaining part of this article we will
most of the time identify tangent vectors and $1$-forms on $M$ using
the isomorphism induced by the metric $g$. 

Let $(g_t,J_t,\o_t,\psp_t,\psm_t)$ be a smooth family of $\SU_3$ structures on
$M$.  We omit the index $t$ when the above tensors are evaluated at
$t=0$, and we use the dot to denote the derivative at $t=0$ in the
direction of $t$.

We start with the study of the $1$-jet at $t=0$ of the family of
$\U_3$ structures $(g_t,J_t,\o_t)$.

\begin{elem} \label{l1}There exist a vector field $\xi$, a section
  $S$ of $\Sym^-M$ and a section $h$ of $\Sym^+M$ $(${\em i.e.} a
  symmetric endomorphism commuting with $J)$, such that
\beq\label{e1}\dot J=JS+\psp_\xi,\eeq
\beq\label{e2}\dot g=g((h+S) \.,\.),\eeq
\beq\label{e3}\dot\o=\f+\xi\i\psp,\eeq
where $\f$ is the $(1,1)$-form defined by $\f(\.,\.)=g(hJ\.,\.)$.
\end{elem}
\bp
Differentiating the identity $J_t^2=-\id_{TM}$ yields that $\dot J$
anti-commutes with $J$. The first formula thus follows directly from
Lemma \ref{ss}.

We write $g_t(\.,\.)=g(f_t \.,\.)$ so $\dot g(\.,\.)=g(\dot
f \.,\.)$. Let $\dot f=h+h'$ be the decomposition of $\dot f$ in
$J$-invariant and $J$-anti-invariant parts. We have to check that $h'=S$.
Differentiating the identity $g_t(\.,\.)=g_t(J_t\.,J_t\.)$ at $t=0$
yields
\bea g(hX,Y)+g(h'X,Y)&=&\dot g(X,Y)=\dot g(JX,JY)+g(\dot JX,JY)+g(JX,\dot JY)\\
&=&g(hJX,JY)+g(h'JX,JY)+g(JSX,JY)\\&&+g(\psp_\xi
X,JY)+g(JX,JSY)+g(JX,\psp_\xi Y)\\
&=&g(hX,Y)-g(h'X, Y)+2g(SX,Y),
\eea
showing that $h'=S$.
The last formula is a direct consequence of (\ref{e1}) and (\ref{e2}):
\bea\dot\o(X,Y)&=&\dot g(JX,Y)+g(\dot
JX,Y)=g((h+S)JX,Y)+g((JS+\psp_\xi)X,Y)\\
&=&g(hJX,Y)+(\xi\i\psp)(X,Y).\eea
\r

This result actually says that the tangent space to the set of all
$\U_3$ structures on $M$ at $(g,J,\o)$ is parametrized by a section
$(\xi,S,\f)$ of the bundle $TM\oplus \Sym^-M\oplus \L^{(1,1)}M$. We now
go forward and describe the $1$-jet of a family of $\SU_3$
structures. Since the reduction from a $\U_3$ structure to a $\SU_3$
structure is given by a section in some $S^1$-bundle, it is not
very surprising that the extra freedom in the tangent space
is measured by a real function ($\mu$ in the notation below):

\begin{elem}\label{l2}
The derivatives at $t=0$ of $\psp_t$ and $\psm_t$ are given by
\beq\label{e4}\dpsp=-\xi\wedge\o+\l\psp+\mu\psm-\frac12 S\bu \psp,\eeq
\beq\label{e5}\dpsm=-J\xi\wedge\o-\mu\psp+\l\psm-\frac12 S\bu \psm,\eeq
where $\l=\frac14 \tr(h)$ and $\mu$ is some smooth function on $M$.
\end{elem}
\bp
By Lemma \ref{ps}, we can write
\beq\label{gt}\dpsp=\a\wedge\o+\l\psp+\mu\psm+Q\bu \psp,\eeq
for some functions $\l,\mu$, $1$-form $\a$ and section $Q$ of $\Sym^-M$.

The fact that $\psp_t$ defines -- in addition to the
$\U_3$ structure $(g_t,J_t)$ -- a $\SU_3$ structure
is characterized by the two equations
\beq\label{su3}g_t(\psp_t,\psp_t)=4\qquad\hbox{and}\qquad
\psp_t(J_tX,Y,Z)=\psp_t(X,J_tY,Z).
\eeq
We consider the symmetric endomorphism $f_t$ introduced above, which
corresponds to $g_t$ in the ground metric $g$. Since the identity acts
on $3$-forms by $-3\id$, the first part of
(\ref{su3}) reads $g(f_t\bu \psp_t,\psp_t)=-12$. Differentiating this at
$t=0$
and using the fact that $\psp$ and $S\bu \psp$ live in orthogonal
components of $\L^3M$, we obtain
\bea 0&=&g(\dot f\bu \psp,\psp)-6g(\dpsp,\psp)=g((h+S)\bu \psp,\psp)-24\l\\
&\stackrel{(\ref{hps})}{=}&6\tr(h)-24\l.\eea
This determines the function $\l$.
We next differentiate the identity $\psp_t\wedge\o_t=0$ at
$t=0$. Since the wedge product with $\o$ vanishes on $\psp$, $\psm$
and on $Q\bu \psp\in\L^{(2,1)+(1,2)}_0M$, we get
$$0=\dpsp\wedge\o+\psp\wedge \dot\o=\a\wedge\o
^2+\psp\wedge(\f+\xi\i\psp)
\stackrel{(\ref{f1})}{=}(\a+\xi)\wedge\o
^2,$$
showing that $\a=-\xi$.
Finally, we differentiate the second part of (\ref{su3}) at $t=0$:
$$\dpsp(JX,Y,Z)+\psp(\dot JX,Y,Z)=\dpsp(X,JY,Z)+\psp(X,\dot JY,Z).$$
Using (\ref{e1}) and (\ref{gt}), this is equivalent to the fact that
the expression
$$-(\xi\wedge\o)(JX,Y,Z)+(Q\bu \psp)(JX,Y,Z)+\psp(JSX,Y,Z)+\psp(\psp_\xi
X,Y,Z)$$
is skew-symmetric in $X$ and $Y$. It is easy to check that
$$-(\xi\wedge\o)(JX,Y,Z)+\psp(\psp_\xi X,Y,Z)=(J\xi\wedge\o)(X,Y,Z),$$
therefore the above condition reduces to
$$(Q\bu \psp)(JX,Y,Z)+\psp(JSX,Y,Z)=(Q\bu \psp)(X,JY,Z)+\psp(X,JSY,Z).$$
Using (\ref{np}), this last relation becomes
\beq\label{sd}\psp((2Q+S)JX,Y,Z)=\psp(X,(2Q+S)JY,Z),\qquad\forall
X,Y,Z\in TM.\eeq
The set of all elements of the form $2Q+S$ satisfying the above
relation is a $\SU_3$-invariant
subspace of $\Sym^-M$. But $\Sym^-M$ is irreducible, and not every
element of 
$\Sym^-M$ satisfies (\ref{sd}) (to see this, just pick up any element in
$\Sym^-M$ and make a direct check). This shows that $2Q+S=0$.

Finally, the relation (\ref{e5}) is a straightforward consequence of
(\ref{e4}). We simply
differentiate the formula $J_t\bu \psp_t=3\psm_t$ (obtained from (\ref{psm}))
at $t=0$ and compute.
\r

Summarizing, we have shown that the tangent space to the set of all
$\SU_3$ structures on $M$ at $(g,J,\o,\psp,\psm)$ is parametrized by a section
$(\xi,S,\f,\mu)$ of the vector bundle $ TM\oplus \Sym^-M\oplus
\L^{(1,1)}M\oplus\RM$.

Let $\a:\L^2M\to TM$ denote the metric adjoint of the linear map
$$X\in TM\mapsto X\i\psp\in\L^2M.$$
A simple check shows that
\beq\label{alp}\a(X\i\psp)=2X,\qquad\a(X\i\psm)=-2JX,\qquad
\a(\tau)=0\ \ \ 
\forall\ \tau\in\L^{(1,1)}M.\eeq
Using the map $\a$, we derive a useful relation between the
components $\dpsp$ and $\dot\o$ of any infinitesimal 
$\SU_3$ deformation:
\beq\label{rr}\L\dpsp\stackrel{(\ref{cont})}{=}-2\xi\stackrel{(\ref{alp})}{=}
-\a(\dot\o).
\eeq

\vs

\section{Deformations of Gray manifolds}

\begin{ede} A {\em Gray structure} on a $6$-dimensional
  manifold $M$ is a $\SU_3$ structure $G:=(g,J,\o,\psp,\psm)$ which
  satisfies the exterior differential system
\beq\label{gray}\begin{cases}d\o=3\psp\\d\psm=-2\o\wedge\o\end{cases}\eeq
A {\em Gray manifold} is a compact manifold endowed with a Gray structure.
\end{ede}

It follows from the work of Reyes-Carri\'on \cite{rc} that a Gray
manifold is a strictly nearly K\"ahler $6$-dimensional compact
manifold with scalar curvature $\scal=30$. We refer to
\cite{gr} for an introduction to nearly K\"ahler geometry. We will use
later on the relations
$$\nabla_X\omega = X\i\psp,\qquad \n_X\psp=-X\wedge\o,$$
which show that $\n\o$ and $\n\psp$ are $\SU_3$-invariant tensor
fields on $M$. Moreover the second equation immediately implies that
$\psi^+$  and $\omega$ are both eigenforms of the Laplace operator for the
eigenvalue $12$. 

Let now $M$ be a compact $6$-dimensional manifold with some Gray
structure $G$ on it. We denote by $\M$ the connected component of
$G$ in the space of Gray structures on $M$. Let $\D$ be the group of
diffeomorphisms of $M$ isotopic to the identity. This group acts on
$\M$ by pull-back and the orbits of this action form the {\em
moduli space of deformations} of $G$.

At the infinitesimal level, the $1$-jet of a curve of Gray
structures $(g_t,J_t,\o_t,\psp_t,\psm_t)$ at $G$ is a tuple
$\gamma:=(\dot g,\dot J,\dot \o,\dpsp,\dpsm)$ determined by a section
$(\xi,S,\f,\mu)$ of the bundle $TM\oplus \Sym^-M\oplus
\L^{(1,1)}M\oplus\RM$ {\em via}
(\ref{e1})-(\ref{e5}), which satisfies the linearized system of (\ref{gray})
\beq\label{dot}
\begin{cases}d\dot\o=3\dpsp,\\
d\dot\psm=-4\dot\o\wedge\o.
\end{cases}\eeq
The space of all tuples $\gamma$ is called the
{\em virtual tangent space} of $\M$ at $G$ and is denoted by $T_G\M$.
The Lie algebra $\chi(M)$ of $\D$ maps to $T_G\M$ by
$X\mapsto \Lie_X G$. Its image, denoted by $\chi ^G (M)$, is a vector
space isomorphic to $\chi(M)/{\mathfrak K}(g)$, where ${\mathfrak
  K}(g)$ denotes the set of Killing vector fields on $M$ with respect
to $g$. The {\em space of infinitesimal Gray deformations} of $G$ is,
by definition, the vector space quotient  $T_G\M/\chi ^G (M)$.

The main purpose of this section is to give a precise description of
this space. 

\begin{ath}\label{infd} Let $G:=(g,J,\o,\psp,\psm)$ be a Gray
  structure on a manifold $M$ such that $(M,g)$ is
not the round sphere $(S^6, can)$. Then the space of infinitesimal
deformations of $G$ 
is isomorphic to the eigenspace for the
eigenvalue $12$ of the restriction of the Laplace operator $\Delta$ to
the space of {\em co-closed} primitive $(1,1)$-forms $\L^{(1,1)}_0M$.
\end{ath}
\bp
A simple but very useful remark is that (except on the round sphere
$S^6$), a Gray structure is completely determined by its underlying Riemannian
metric. The reason for that is the fact that the metric defines
(locally) a unique line of Killing spinors with positive Killing
constant, which, in turn, defines the almost complex structure, and,
together with the exterior derivative of the K\"ahler form, one
recovers the whole $\SU_3$ structure.

By the Ebin Slice Theorem \cite{be}, each infinitesimal
deformation of $G$ has a unique representative $\gamma=(\dot g,\dot J,\dot
\o,\dpsp,\dpsm)\in T_G\M$ such that 
\beq\label{ebin}\d\dot g=0\qquad\hbox{and}\qquad \tr_g\dot g=0.\eeq
Let $(\xi,S,\f,\mu)$ be the section of the bundle $TM\oplus \Sym^-M\oplus
\L^{(1,1)}M\oplus\RM$ determined by $\gamma$ {\em via} the equations
(\ref{e1})-(\ref{e5}). We have to interpret the system
(\ref{dot})-(\ref{ebin}) in terms of $(\xi,S,\f,\mu)$. The main idea
is to show that $\mu$ and $J\xi$ are eigenvectors of the Laplace
operator on functions and $1$-forms corresponding to eigenvalues
which do not belong to its spectrum, thus forcing them to vanish.

We start by taking the exterior product with $\psp$ in the first equation of
(\ref{dot}) and use (\ref{e4}) to get
$$d\dot\o\wedge\psp=3\dpsp\wedge\psp=3\mu\psm\wedge\psp=-12\mu dv.$$
On the other hand, using (\ref{e3}) and taking (\ref{f3}) into account
yields
$$d\dot\o\wedge\psp=d(\f+\xi\i\psp)\wedge\psp=d((\f+\xi\i\psp)\wedge\psp)
=d(\xi\wedge\o ^2),$$
whence
\beq\label{mu} -12\mu dv=d(\xi\wedge\o ^2).\eeq
We apply the contraction $\L$ to the first equation of
(\ref{dot}) and use (\ref{cont}), (\ref{e3}) and (\ref{e4}):
\beq\label{inter}
-6\xi=3\L\dpsp=\L d\dot\o=\L d\f+\L d(\xi\i\psp).\eeq
In order to compute the last term, we apply the general formula
(\ref{rr}) to the particular deformation of the $\SU_3$ structure
defined by the flow of $\xi$:
$$\L(\Lie_\xi\psp)=-\a(\Lie_\xi\o).$$
Since $d\o=3\psp$ and $d\psp=0$, we get
$$\L d(\xi\i\psp)=\L(\Lie_\xi\psp)=-\a(\Lie_\xi\o)
=-\a(d(\xi\i\o))-\a(\xi\i
d\o)\stackrel{(\ref{alp})}{=}-\a(dJ\xi)-6\xi,$$
which, together with (\ref{inter}), yields
\beq\label{eq2}\L(d\f)=\a(dJ\xi).\eeq
We now examine the second equation of the system (\ref{dot}). From
(\ref{e5}) we get
$$d\dpsm=-dJ\xi\wedge\o+3
J\xi\wedge\psp-d\mu\wedge\psp+d\l\wedge\psm-2\l\o ^2-\frac12
d(S\bu \psm).$$
We apply the contraction $\L$ to this formula and use the second
equation of (\ref{dot}) together 
with  (\ref{lam}) and (\ref{lam1}):
\beaa\label{aa} -4\dot\o&=&\L(-4\dot\o\wedge\o)+4(\L\dot
\o)\o=\L d\dpsm+4(\L\dot \o)\o\\ 
&=&-dJ\xi-\L(dJ\xi)\o-3\xi\i\psp-Jd\mu\i\psp+Jd\l\i\psm\\&&-8\l\o-\frac12
\L d(S\bu \psm)+4(\L\dot \o)\o.
\eeaa
Applying $\a$ to this equality and using (\ref{psm}), (\ref{e3}) and
(\ref{rr}) yields
$$-8\xi=-\a(dJ\xi)-6\xi-2Jd\mu+2d\l-\frac12
\a\L d(S\bu \psm).
$$
From (\ref{eq2}) we then get
\beq\label{k1}\xi=Jd\mu+\frac12\L d\f-d\l+\frac14
\a\L d(S\bu \psm).\eeq
\begin{elem}\label{cl} Let $(g,J,\o,\psp,\psm)$ be a Gray structure on
a manifold $M$. For every section
$S$ of $\Sym^-M$ and $(1,1)$-form $\f$, the following relations hold:
\beq\label{claim}\L d\f=\d h+2d\l\eeq
\beq\label{claim3}\d h=-J\d\f\eeq
\beq\label{claim1}\d(S\bu \psp)=-\L d(S\bu \psm)-2\d S\i\psp\eeq
\beq \label{claim2}\a\L d(S\bu \psm)=-2\d S \eeq
\beq\label{claim4}\L\d(S\bu \psp)=0\eeq
where $h$ is the endomorphism defined by $\f(\.,\.)=g(hJ\.,\.)$ and
$\l=\frac14 \tr (h)=\frac12 \L\f.$ In the above formulas, $\d$ stands
for the usual co-differential when applied to an exterior form, and for
the divergence operator $($cf. \cite{ab}, $1.59)$ when acting on
symmetric tensors. 
\end{elem}
\ni Since the proof is rather technical, we postpone it to the end of this
section.

\ni Using (\ref{e2}), (\ref{claim}), (\ref{claim3}) and
(\ref{claim2}), the relation 
(\ref{k1}) becomes 
\beq\label{xi1}\xi= Jd\mu+\frac12\d h-\frac12\d S=
Jd\mu-J\d\f-\frac12\d \dot g.
\eeq
From (\ref{ebin}) and (\ref{xi1}) we obtain
\beq\label{xi}\xi= Jd\mu-J\d\f.
\eeq
On the other hand, (\ref{xi})
shows that $\mu$ is an eigenfunction of $\Delta$ for the eigenvalue 6:
\beq\label{xi2}\Delta\mu=\d d\mu=\d(J\xi)=-*d*(J\xi)=-\frac12*d(\xi\wedge\o ^2)
\stackrel{(\ref{mu})}{=}6\mu.\eeq
Now, the Obata theorem
(cf. \cite[Theorem 3]{ob}) says that on a compact $n$-dimensional
Einstein manifold of positive scalar curvature $\scal$,
every eigenvalue of the Laplace operator is greater than
or equal to $\scal/(n-1)$, and equality can only occur on the standard sphere.
Since $(M,g)$ is not isometric to the standard
sphere and is Einstein with scalar curvature $\scal=30$, the Obata theorem
thus implies that $\mu=0$. Since $\l=\frac14\tr(h)=\frac14\tr(\dot g)$, the
second part of (\ref{ebin}) also shows that $\l=0$.  
Taking (\ref{ebin}), (\ref{claim3}), (\ref{claim1}) and (\ref{xi})
        into account, the equation (\ref{aa}) now becomes 
\bea -4\f-4\xi\i\psp&=&-4\dot\o=-dJ\xi-\L(dJ\xi)\o-3\xi\i\psp-\frac12
\L d(S\bu \psm)\\
&=&-dJ\xi-\L(dJ\xi)\o-
3\xi\i\psp+\frac12\d(S\bu \psp)-\xi\i\psp,  
\eea
whence 
$$dJ\xi=-\L(dJ\xi)\o+4\f+\frac12\d(S\bu \psp).$$
Applying $\L$ to this relation and using (\ref{claim4}) yields
$\L(dJ\xi)=3\L(dJ\xi)$, so $dJ\xi=4\f+\frac12\d(S\bu \psp).$
By (\ref{xi}) we have $J\xi=\d\f$, so
$$\Delta(J\xi)=\d dJ\xi=4\d\f=4J\xi,$$
{\em i.e.} $J\xi$ is an eigenform of the Laplace operator with
eigenvalue 4.
On the other hand, the Bochner formula on $1$-forms, 
$$\Delta=\n^*\n+\Ric=\n^*\n+5\id$$ 
shows, by integration over $M$, that 4 cannot be an eigenvalue of
$\Delta$ , so
$\xi$ has to vanish identically.
From Lemmas \ref{l1} and \ref{l2}, we get
$$
\dot \psi^+  =  - \frac12 S \bu  \psp,\qquad
\dot \psi^-  = - \frac12 S \bu  \psm,\qquad
\dot \omega  =  \f \in \Lambda^{(1,1)}_0M \ .
$$
Plugging these equations into (\ref{dot}) yields
\beq\label{dot2}
\begin{cases}d\f=-\frac32 S\bu \psp,\\
\d(S\bu \psp)=-*d*(S\bu \psp)\stackrel{(\ref{sp})}{=}*d(S\bu
\psm)=-2*d\dpsm
\stackrel{(\ref{dot})}{=}8*(\dot\o\wedge\o)\stackrel{(\ref{sp1})}{=}-8\f.
\end{cases}\eeq
Thus $\f$ is a co-closed eigenform of the Laplace
operator for the eigenvalue 12.

Conversely let us assume that $\f \in \Omega^{(1,1)}_0M$
is co-closed and satisfies $\Delta \f = 12\f$. We have
to show that $\f$ defines an infinitesimal deformation of
the given Gray structure. The main point is to remark that $d\f$
is a form in $\Lambda^{(2,1)+(1,2)}_0M$. 

\begin{elem} If $\f$ is a co-closed form in $\O ^{(1,1)}_0M$ then
$d\f\in \O ^{(2,1)+(1,2)}_0M.$
\end{elem}
\bp
Using Lemma \ref{ps}, this amounts to check
that $d\f$ satisfies the system 
\beq
\begin{cases}d\f\wedge\psp=0,\\
d\f\wedge\psm=0,\\
\la d\f, X \wedge \omega\ra =0,\qquad\forall\ X\in TM.
\end{cases}\eeq
Each of these equations follows easily:
$$d\f\wedge\psp=d(\f\wedge\psp)=0,$$
$$d\f\wedge\psm=d(\f\wedge\psm)-\f\wedge
d\psm\stackrel{(\ref{gray})}{=}2\f\wedge\o
^2\stackrel{(\ref{star})}{=}4\la\f,\o\ra=0,$$ 
$$\la d\f, X \wedge \omega\ra =\la \L d\f,
X\ra\stackrel{(\ref{claim})}{=}-\la J\d\f,X\ra dv=0.$$
\r

Consequently, there exists a
unique section $S$ of $\Sym^-M$ with $d\f = -\frac32 S\bu \psp$. Taking
$\xi=0$, $\l=0$ and $\mu=0$, the equations (\ref{e1})-(\ref{e5})
define an infinitesimal $\SU_3$-deformation by $\dot \omega :=
\varphi$, $\dot \psi^+: = -\frac12S\bu \psp$ 
and $\dot \psi^-: = -\frac12S\bu \psm$. It remains to show that
$(\dot\o,\dpsp,\dpsm)$ satisfy the linearized system~(\ref{dot}).
The first equation is clear by definition and the second
is equivalent to $d (S\bu \psm) = 8\f \wedge \omega$. Using again
(\ref{sp}) and (\ref{sp1}), this last equation is equivalent to
$\d(S\bu \psp)=-8\f$, which follows directly from the hypothesis on $\f$
and the definition of $S$.
\r
\vs

\ni{\em Proof of Lemma \ref{cl}}. 
By (\ref{ii}), $\n_{X}J=\psp_X\in \L^{(0,2)+(2,0)}M$ so
$\la\f,\n_{X}J\ra=0$ for all vectors $X$. Identifying $\f$ with the
corresponding skew-symmetric endomorphism of $TM$, we compute in a
local orthonormal basis $\{e_i\}$ parallel at some point:
\bea \L d\f&=&\frac12 Je_i\i e_i\i(e_k\wedge\n_{e_k}\f)=\frac12
Je_i\i\n_{e_i}\f-\frac12 
Je_i\i(e_k\wedge(\n_{e_k}\f)e_i)\\
&=&\frac12(\n_{e_i}\f)Je_i-\frac12(\n_{Je_i}\f)e_i+\frac12
e_k(\n_{e_k}\f)(e_i,Je_i)\\
&=&(\n_{e_i}\f)Je_i+d\la\f,\o\ra-e_k\f(e_i,(\n_{e_k}J)e_i)\\
&=&(\n_{e_i}(\f J))e_i-\f((\n_{e_i}J)e_i)+2d\l-2e_k\la\f,\n_{e_k}J\ra\\
&=&-(\n_{e_i}h)(e_i)+2d\l=\d h+2d\l. \eea
This proves (\ref{claim}). To prove the second relation, we notice that 
$$g((\n_{e_i}J)(\f e_i),X)=\psp(e_i,\f e_i,X)=2\la\psp_X,\f\ra=0,$$
so
$$\d h=-(\n_{e_i}h)e_i=(\n_{e_i}(J\f))e_i=
-J\d\f+(\n_{e_i}J)\f e_i=-J\d\f.$$
We will now use several
times (and signal this by a star above the equality sign in the
calculations below) the fact that any $\SU_3$-equivariant map
$\Sym^-M\to\L^2M$ is automatically zero, by the Schur Lemma. For every
$3$-form $\tau$ we can write
\bea \L d(\tau)&=&\frac12 Je_i\i e_i\i(e_k\wedge\n_{e_k}\tau)=
e_k\wedge(\frac12 Je_i\i e_i\i\n_{e_k}\tau)+Je_k\i\n_{e_k}\tau\\
&=&d\L\tau-\frac12 e_k\wedge(\n_{e_k}J)e_i\i e_i\i\tau
+Je_k\i\n_{e_k}\tau. \eea
In particular, for $\tau=S\bu \psp$ we have $\L\tau=0$ and
$e_k\wedge(\n_{e_k}J)e_i\i e_i\i\tau=0$, so from (\ref{a}) and the
remark above we get
\bea \L d(S\bu \psm)&=&Je_k\i\n_{e_k}(S\bu \psm)= 
Je_k\i((\n_{e_k}S)\bu \psm)+Je_k\i(S\bu (\n_{e_k}\psm))\\
&\stackrel{*}{=}&-Je_k\i((\n_{e_k}S)e_j\wedge e_j\i\psm)\\
&=&-\la (\n_{e_k}S)e_j,Je_k\ra e_j\i\psm+(\n_{e_k}S)e_j\wedge
(Je_k\i e_j\i\psm)\\
&\stackrel{*}{=}&\la (\n_{e_k}S)Je_j,e_k\ra Je_j\i\psp+(\n_{e_k}S)e_j\wedge
(Je_k\i e_j\i\psm)\\
&=&-\d S\i\psp+(\n_{e_k}S)e_j\wedge (Je_k\i
e_j\i\psm), \eea
so by (\ref{psm})
\beq\label{fi}\L d(S\bu \psm)=-\d S\i\psp+(\n_{e_k}S)e_j\wedge (e_k\i
e_j\i\psp).
\eeq
On the other hand, we have 
\bea\d(S\bu \psp)&=&e_k\i\n_{e_k}(Se_j\wedge e_j\i\psp)\stackrel{*}{=}
e_k\i((\n_{e_k}S)e_j\wedge e_j\i\psp)\\
&=&-\d S\i\psp-(\n_{e_k}S)e_j\wedge e_k\i e_j\i\psp\\
&\stackrel{(\ref{fi})}{=}&-2\d S\i\psp-\L d(S\bu \psm),\eea
thus proving (\ref{claim1}).
Let $\sigma:=(\n_{e_k}S)e_j\wedge (e_k\i e_j\i\psp)$ denote the
last summand in (\ref{fi}). A similar calculation easily shows that $J\bu
\s=0$, so $\s$ belongs to $\L^{(1,1)}M$. From (\ref{alp}) we thus get
$$\a\L d(S\bu \psm)=-\a(\d S\i\psp)=-2\d S.$$
Finally, the relation (\ref{claim4}) can be checked in the same way:
\bea \L\d(S\bu \psp)&=&-\frac12 Je_k\i e_k\i(e_j\i \n_{e_j}(S\bu \psp))\\
&=&-\frac12\bigg(e_j\i \n_{e_j}(Je_k\i e_k\i(S\bu \psp))-e_j\i
(\n_{e_j}J)e_k\i e_k\i(S\bu \psp))\bigg)\\
&\stackrel{*}{=}&\d\L(S\bu \psp)=0.
\eea
This proves the lemma.
\hfill $\square$

\section{Concluding Remarks}

So far we have identified the space of infinitesimal deformations
of a given Gray structure with the space of co-closed primitive $(1,1)$-forms
which are eigenforms of the Laplace operator for the eigenvalue
12. In order to proceed further there are two immediate
options. One could try to
compute the second  derivative of a curve of Gray structures
and obtain additional equations. However, this leads to
quadratic expressions, which for the moment seem to be rather difficult 
to handle.

A second natural task is to consider the known homogeneous examples
and to study 
the question whether or not there exist at least infinitesimal
deformations. This amounts to study the Laplace operator
on 2-forms on certain homogeneous spaces, which should reduce to a tractable
algebraic problem.

Since a nearly K\"ahler deformation always gives rise to an Einstein
deformation, we could equally well ask for the existence of 
infinitesimal Einstein deformations on the known nearly K\"ahler
manifolds. Except for the standard sphere, where no deformations 
exist, there seems to be nothing known in this direction.

The situation for the standard sphere is quite different. 
Th. Friedrich has shown in \cite{fr} that the action of the isometry
group $\SO_7$ on the set of nearly K\"ahler
structures on the round sphere $(S^6,can)$ is transitive. The
isotropy group of this action at the standard nearly 
K\"ahler structure is easily seen to be the group $G_2$ (the
stabilizer in $\SO_7$ of a vector cross product). The space of nearly K\"ahler
structures  
on the round sphere is thus isomorphic to $\SO_7/G_2\simeq
\RP^7$ (see \cite[Prop. 7.2]{jbb}). More geometrically, the set of nearly K\"ahler
structures compatible with the round metric on $S^6$ can be identified
with the set of non-zero real Killing spinors (modulo constant
rescalings), and 
it is well-known that the space of real Killing spinors on $(S^6,can)$ is
isomorphic to $\RM^8$.

The counterpart of Theorem \ref{infd} on the round sphere $(S^6,can)$
can be stated as follows: 

\begin{ath}\label{infd1} Let $G:=(g,J,\o,\psp,\psm)$ be a Gray
structure on $S^6$ such that $g$ is the round metric. Then the space
of infinitesimal 
deformations of $G$ 
is isomorphic to the eigenspace for the
eigenvalue $6$ of the Laplace operator $\Delta$ on functions, and is,
in particular, $7$-dimensional.
\end{ath}
\proof
Since there are no Einstein deformations on the sphere, we may assume
$\dot g = 0$ from the beginning. From (\ref{e2}) we get $\f = 0, h=0, S = 0$
and in particular $\lambda = 0$ too. Then equation~(\ref{xi})
gives $\xi = Jd\mu$ and (\ref{xi2}) shows that $\mu$ is an
eigenfunction of $\Delta$ on $(S^6, can)$
corresponding to its first non-zero eigenvalue, $6$. These
eigenfunctions (also called {\em first 
spherical harmonics}) satisfy $\nabla_X d\mu
=-\mu X $ for all tangent vectors $X$.  We define the infinitesimal
$SU_3$ deformation $\dot \omega  := \xi 
\lrcorner \psi^+$,  
$\dot \psi^+ := -\xi \wedge \omega +\mu \psi^-$ and $\dot \psi^- :=
-J\xi \wedge \omega -\mu \psi^+$. A short calculation easily shows
that this indeed is a solution of the linearized system~(\ref{dot}).
\qed

\labelsep .5cm

\end{document}